\documentclass[letterpaper,11pt,openany,oneside]{article}
\usepackage[T1]{fontenc}
\usepackage[latin1]{inputenc}
\usepackage{epsfig}
\usepackage{amsmath,amssymb}
\usepackage[active]{srcltx}
\usepackage{dsfont}
\usepackage[titles]{tocloft}
\usepackage{amssymb}
\usepackage{tikz}
\usetikzlibrary{decorations.pathreplacing}
\numberwithin{equation}{section}
\newtheorem{theoreme}{Theorem}[section]

\newtheorem{claim}[theoreme]{Claim}

\newtheorem{lemme}{Lemma}
\newtheorem{thm}{Theorem}[section]
\newtheorem{cor}[thm]{Corollary}
\newtheorem{lem}[thm]{Lemma}
\newtheorem{pro}[thm]{Proposition}
\newtheorem{defi}[thm]{Definition}
\newtheorem{rem}[thm]{Remark}

\newtheorem{exa}[thm]{Example}

\newtheorem{countexa}[thm]{Counter-example}

\newenvironment{proof}[1][Proof]{\noindent \textbf{#1.}~ }
{\hfill\rule{2mm}{2mm} \vspace{\parskip} }

\newcommand{\RR}{\ensuremath{\mathbb R}}

\newcommand{\NN}{\ensuremath{\mathbb N}}

\newcommand{\N}{\ensuremath{\mathcal N}}

\newcommand{\T}{\ensuremath{\mathcal T}}

\newcommand{\deriv}{\mathrm{d}}

\newcommand{\De}{\Delta}

\title{Limit value for optimal control with general means\thanks{This research was partially supported by Agence National de la Researche (grant ANR-10-BLAN 0112). This article is done as part of the PhD thesis of the first author, he wishes to thank his supervisor Sylvain Sorin for advising and comments.}}
\author{Xiaoxi LI \thanks{CNRS, IMJ-PRG, UMR 7586, Sorbonne Universit\'es, UPMC Univ Paris 06, Univ Paris Diderot, Sorbonne Paris Cit\'e, Case 247, 4 Place Jussieu, 75252 Paris, France.},
\  Marc QUINCAMPOIX \thanks{Laboratoire de Math\'ematiques de Bretagne Atalantique, UMR 6205, Universit\'e de Brest, 6 Avenue Le Gorgeu, 29200 Brest, France.},
\ J\'er\^ome RENAULT \thanks{TSE (GREMAQ, Universit\'e Toulouse 1 Capitole), 21 all\'ee de Brienne, 31000 Toulouse, France.}}

\date{\today }

\begin{document}
\maketitle

\begin{abstract}
We consider optimal control problem with an integral cost which is a mean of a given function. As a particular case, the cost concerned is the Ces\`aro average. The limit of the value with Ces\`aro mean when the horizon tends to infinity is widely studied in the literature. We address the more general question of the existence of a limit when the averaging parameter converges, for values defined with means of general types.

We consider a given function and a family of costs defined as the mean of the function with respect to a family of probability measures -- the evaluations -- on $\RR_+$. We provide conditions on the evaluations in order to obtain the uniform convergence of the associated value function (when the parameter of the family converges).

Our main result gives a necessary and sufficient condition in term of the total variation of the family of probability measures on $\RR_+$. As a byproduct, we obtain the existence of a limit value (for general means) for control systems having a compact invariant set and satisfying suitable nonexpansive property.

\end{abstract}


\newpage

\section{Introduction}
We consider a control system defined on $\RR^d$ whose dynamic is given by
\begin{align} \label{eq:dynamics}
y'(t)=f\big(y(t),u(t)\big)
\end{align}
where $f:\RR^d\times U \to \RR^d$ and $u(\cdot)$ is a measurable function -- called the control -- from $\RR_+$ to $U$ a fixed metric space. We will make later on assumptions on ($\ref{eq:dynamics}$) ensuring that for any initial condition $y(0)=y_0$, and any measurable control $u(\cdot)$, the equation ($\ref{eq:dynamics}$) has a unique solution $t\mapsto y(t,u,y_0)$ defined on $\RR_+$.

To any pair \big($y_0,u(\cdot)$\big), we associate a cost
\begin{align*}
\int_{0}^{+\infty}g\big(y(t,u,y_0),u(t)\big)\deriv\theta(t),
\end{align*}
where $g:\RR^d\times U\to \RR$ is Borel measurable bounded and $\theta\in\De(\RR_+)$ is a Borel probability measure on $\RR_+$. We called $\theta$ an \textit{evaluation} throughout the article.

We will refer to the previously described optimal control problem by the short notation $\mathcal{J}=\langle U,g,f\rangle$. Let $\theta\in\De(\RR_+)$, we define for $\mathcal{J}=\langle U,g,f\rangle$ the following value function:
\begin{align}\label{eq:V-theta}
V_\theta(y_0)=\inf_{u(\cdot)\in\mathcal{U}}\int_{0}^{+\infty}g\big(y(t,u,y_0),u(t)\big)\deriv\theta(t),
\end{align}
where $\mathcal{U}$ denotes the set of measurable controls $u:[0,+\infty)\to U$.

\noindent Typical means in the definition ($\ref{eq:V-theta}$) of the value function are well studied in the literature for

\noindent \textit{Ces\`aro  mean}: \ $\forall t>0$, $\overline{\theta}_t$ with density $s\mapsto f_{\overline{\theta}_t}(s)=\frac{1}{t}\mathds{1}_{[0,t]}(s)$, and the $t$-horizon value is
\begin{equation*}
V_{\overline{\theta}_t}(y_0)=\inf_{u(\cdot)\in \mathcal{U}} \frac{1}{t}\int_{0}^{t} g\big(y(s,u,y_0),u(s)\big) \deriv s
\end{equation*}
\noindent \textit{Abel mean}: \ $\forall \lambda\in(0,1]$, $\theta_\lambda$ with density $s\mapsto f_{\theta_\lambda}(s)=\lambda e^{-\lambda s}$, and the $\lambda$-discounted value is
\begin{equation*}
V_{\theta_\lambda}(y_0)=\inf_{u(\cdot)\in \mathcal{U}}\int_{0}^{+\infty} \lambda e^{-\lambda s}g\big(y(s,u,y_0),u(s)\big) \deriv s
\end{equation*}

The limit of the above value functions as $t$ tends to infinity or as $\lambda$ tends to zero are well investigated in the control literature, (cf. \cite{Alvarez_10}, \cite{Arisawa_98a}, \cite{Arisawa_98b}, \cite{Bensoussan_88}, \cite{Gaitsgory_86}, \cite{Khasminskii_68} and the references therein), which are often called ergodic control.

When $\theta\in\De(\RR_+)$ is given, the contribution of the interval $[T,+\infty)$ in the mean ($\ref{eq:V-theta}$) is less and less significant as $T$ becomes large. Thus the control problem is essentially interesting only on $[0,T_0]$ for certain $T_0$, which we roughly name the "duration" for the problem. In this article, we are interested in the long-run property of $\mathcal{J}$, $i.e.$, the asymptotic behavior of the function $\theta\mapsto V_\theta$ when the "duration" of $\theta$ tends to infinity.  In the particular examples of Ces\`aro mean and Abel mean, the uniform convergence of $V_{\overline{\theta}_t}$ as $t$ tends to infinity and of $V_{\theta_\lambda}$ as $\lambda$ tends to $0$ are studied.  It is a priori unclear how to define the "duration" of a general evaluation $\theta$ over $\RR_+$. If one just assumes the expectation of $\theta$ to be large, we can obtain very different value functions, as is shown by the following

\begin{exa}\label{exa:0-1} Consider the uncontrolled dynamic $y(t)=t$, the running cost $t\mapsto g(t)=\mathds{1}_{\cup_{m=1}^{\infty}[2m-1, 2m]}(t)$, and two sequences of evaluations $(\mu^k)_{k\geq 1}$ and $(\nu^k)_{k\geq 1}$ with densities: $f_{\mu^k}=\frac{1}{k}\mathds{1}_{\cup_{m=1}^k[2m-1, 2m]}$ and $f_{\nu^k}=\frac{1}{k}\mathds{1}_{\cup_{m=1}^k[2m-2, 2m-1]}$. Clearly, $V_{\mu^k}=1$ and $V_{\nu^k}=0$, $\forall k\geq 1$.
\end{exa}

For this reason, we introduce an asymptotic regularity condition for evaluations, called the \textit{long-term condition} (LTC for short), to express the "large duration" and the "asymptotic uniformity of distributions over $\RR_+$", and we will study the convergence of the value functions along a sequence of evaluations satisfying the LTC.

More precisely,  for any $s\geq 0$, we define the \textit{$s$-total variation} of an evaluation $\theta$ to be the total variation between the measure $\theta$ and its $s$-shift along $\RR_+$: $$TV_s(\theta)=\max_{Q\in\mathcal{B}(\RR_+)}|\theta(Q)-\theta(Q+s)|.$$
We say that a sequence of evaluations $(\theta^k)_{k\geq 1}$ satisfies the LTC if:
$$\forall S>0, \ \sup_{0\leq s\leq S}TV_s(\theta^k) \xrightarrow[k\to\infty]\ {} 0.$$
The optimal control problem $\mathcal{J}=\langle U,g,f\rangle$ has a \textit{general limit value} given by some function $V^*$ defined on $\RR^d$ if for any sequence $(\theta^k)_k$ satisfying the LTC, $\left(V_{\theta^k}(y_0)\right)_k$ converges uniformly to  $V^*$ as $k$ tends to infinity.

Our main result (Theorem \ref{thm:main}) states that for any $(\theta^k)_k$ satisfying the LTC, $(V_{\theta^k})_k$ converges uniformly if and only if the family $\{V_{\theta^k}\}$ is totally bounded with respect to the uniform norm. Moreover, in this case, the limit is characterized by the following:
\begin{align} \label{eq:V*}
V^*(y_0)=_{def}\sup_{\theta\in\De({\RR_+})}\inf_{s\in\RR_+} \inf_{u\in\cal{U}}\int_{0}^\infty g\Big(y(t+s, u,y_0),u(t+s)\Big)\deriv\theta(t), \ \  \forall y_0\in\RR^d.
\end{align}
The above function $V^*$ naturally appears to be the unique possible long-term value function of the control problem.

As a byproduct of our main result, we obtain the existence of the general limit value for any control problem $\mathcal{J}=\langle U,g,f\rangle$ with a running cost $g$ that does not depend on $u$ and with a control dynamic ($\ref{eq:dynamics}$) which is non-expansive and has a compact invariant set. This can be viewed as a generalization of already obtained results in \cite{Quincampoix_2011} for optimal control with Ces\`aro mean.

Existing results in the erdogic control literature are concerned mainly with the convergence of the $t$-horizon Ces\`aro mean values or the convergence of the $\lambda$-discounted Abel mean values. To the best of the authors' knowledge, this paper is the first to consider general long-term evaluations for optimal control problems.

Also it is worth pointing out that while many works (including \cite{Alvarez_10},  \cite{Arisawa_98a}, \cite{Arisawa_98b}, \cite{Bensoussan_88}, \cite{Gaitsgory_86}, \cite{Khasminskii_68}) suppose controllability or ergodicity conditions, the present approach does not reply on such conditions. This could be understood by the fact that the limit value $V^*$ may depend on the initial state $y_0$ (which does not occur under ergodic or controlability assumptions).

We also make here a link to the discrete time framework, in which an evaluation $\theta=(\theta_m)_{m\geq 1}$ is a probability measure over positive integers $\NN^*=\NN \backslash\{0\}$, and $\theta_t$ is the weight for the stage-$t$ payoff. The analogue notion of \textit{total variation} is defined for any $\theta\in\De(\NN^*)$: $TV(\theta)=\sum_{m=1}^\infty|\theta_{m+1}-\theta_m|$ (cf. \cite{Sorin_2002} and \cite{Renault_2014}). Recently, the existence of the general limit value of dynamic optimization problems in several discrete time frameworks has been obtained in \cite{Renault_2014}, \cite{Renault_2013} and \cite{Ziliotto_2014}. Our work is partially inspired by \cite{Renault_2014} and similar tool within the proof appeared in \cite{Renault_2011}.

The article is organized as follows. Section \ref{sec:preliminary} contains some preliminary notations and basic examples. The long-term condition is introduced and studied in Section \ref{sec:LTC}. Section \ref{sec:result} contains our main result and its consequences. We discuss in the end of this section two (counter)examples. Section \ref{sec:proof} is devoted to the proof of the main result. A weaker notation of LTC is discussed in Section \ref{sec:LTC'}.

\section{Preliminaries}\label{sec:preliminary}
Consider the optimal control problem $\mathcal{J}=\langle U, g,f\rangle$ described by (\ref{eq:dynamics})-(\ref{eq:V-theta}).
We make the following assumptions on $g$ and $f$:
\begin{equation}\label{ex:assumptions}
\begin{cases}
\text{ the function} \ g:\RR^d\times U\to\RR \text{ is Borel measurable and bounded}; \\
\text{ the function}  \ f:\RR^d\times U\to\RR^d \text{ is Borel measurable, and satisfies:}\\
(*).\ \exists L\geq 0, \forall (y,\overline{y})\in\RR^{2d}, \forall u\in U, ||f(y,u)-f(\overline{y},u)||\leq L||y-\overline{y}||,\\
(**).\ \exists a>0, \forall (y,u)\in\RR^d\times U, ||f(y,u)||\leq a(1+||y||).
\end{cases}
\end{equation}
\noindent Under these hypotheses, given any control $u(\cdot)$ in $\mathcal{U}$ and any initial starting state $y_0\in\RR^d$, (\ref{eq:dynamics}) has a unique absolutely continuous solution $t\mapsto y(t, u, y_0)$ defined on $[0,+\infty)$. As the running cost function $g:\RR^d\times U\to\RR$ is bounded, we can always assume that $g:\RR^d\times U\to [0,1]$ after some affine transformation.

Below we introduce several notations.

\noindent \textbf{$\theta$-evaluated cost $\gamma_\theta(y_0,u)$}
\ \ Given $\theta\in\De(\RR_+)$ and $y_0\in\RR^d$, the $\theta$-evaluated cost induced by a control $u(\cdot)\in \mathcal{U}$ is denoted by:
\begin{align*}
\gamma_\theta(y_0,u)=\int_{[0,+\infty)} g\Big(y(s,u,y_0),u(s)\Big)\deriv\theta(s),
\end{align*}
With this notation, the $\theta$-value function in (\ref{eq:V-theta}) writes as $V_\theta(y_0)=\inf_{u(\cdot)\in \mathcal{U}}\gamma_{\theta}(y_0,u)$.

\noindent \textbf{Reachable map $R_t$}  \ \ For any $y_0\in\RR^d$, the reachable map in $\RR_+$, $t\mapsto R_t(y_0)$, is defined as:
\begin{align}\label{eq:Reach}
R_t(y_0)=\Big\{\overline{y}\in\RR^d\big|\exists \ u(\cdot)\in\mathcal{U}: y(t,u,y_0)=\overline{y}\Big\}.
\end{align}
$R_t(y_0)$ represents the set of states that the dynamic can reach via certain control at time $t$, starting from the initial state $y_0$ at time $0$. We write $R^t(y_0)=\cup_{s=0}^t R_s(y_0)$ and $R(y_0)=\cup_{s=0}^\infty R_s(y_0)$. $R(y_0)$ is the set of states that can be reached at any finite time starting from $y_0$.

\noindent \textbf{Image measure $\T_t\sharp\theta$ and the auxiliary value function $V_{\T_t\sharp\theta}$}\ \ \  Given $t\in\RR$ and $\theta$ in $\De(\RR_+)$, we use $\T_t\sharp\theta$ to denote the image (push-forward) measure of $\theta$ by the function $\T_t: s\mapsto s+t$, $i.e.$,
$$\T_t\sharp\theta (Q)= \theta\left(\T_t^{-1}(Q)\right)=\theta \big((Q-t)\cap \RR_+\big), \ \forall Q\in {\cal B}(\RR_+),$$
where ${\cal B}(\RR_+)$ denotes the set of all Borel subsets in $\RR_+$. This leads us to write the $t$-shift $\theta$-evaluated cost induced by a control $u$ as follow:
\begin{align}\label{eq:gamma+t+shift}
\gamma_{\T_t\sharp\theta}(y_0, u)=\int_{[0,+\infty)} g\Big(y(s+t,u,y_0),u(s+t)\Big)\deriv\theta(s), \ \forall t\geq 0.
\end{align}
Taking on both sides of (\ref{eq:gamma+t+shift}) the infimum over $u(\cdot)\in\mathcal{U}$ and using the notation of reachable map $R_t$, we obtain the $t$-shift $\theta$-value function
 \begin{align}\label{eq:V+t+shift}
 V_{\T_t\sharp\theta } (y_0)=\inf_{u(\cdot)\in \mathcal{U}} \int_{[0,+\infty)} g\Big(y(s+t,u,y_0),u(s+t)\Big) \deriv \theta(s)=\inf_{\overline{y}\in R_t(y_0)}V_{\theta}(\overline{y}).
 \end{align}

\noindent \textbf{s-total variation} Given an evaluation $\theta$, define its $s$-total variation for each $s\geq 0$:
\begin{align} \label{eq:TV_s}
TV_s(\theta)=\sup_{Q\in {\cal B}(\RR_+)} |\theta(Q)-\theta(Q+s)|.
\end{align}

\noindent \textbf{Long-term condition (LTC)} \ A sequence of evaluations $(\theta^k)_{k\geq 1}$ satisfies the LTC if:
\begin{align} \label{eq:LTC}
\forall S>0, \ \overline{TV}_S(\theta^k)=_{def}\sup_{0\leq s\leq S}TV_s(\theta^k)\xrightarrow[k \to \infty]{} \ 0.
\end{align}
In this article, we are concerned with the following notation of limit value for optimal control problems with general means.
\begin{defi}\label{def:glv} Let $V$ be a function defined on $\RR^d$. The optimal control problem $\mathcal{J}$ admits $V$ as the \textbf{general limit value} if: for any sequence of evaluations $(\theta^k)_{k\geq 1}$ satisfying the LTC, for all $y_0$ in $\RR^d$, $\left(V_{\theta^k}(y_0)\right)$ converges to $V(y_0)$ as $k$ tends to infinity, and moreover the convergence is uniform in $y_0$.
\end{defi}

Below are some basic examples of optimal control problems in which the general limit value exists.

\begin{exa} \label{exampe:1} $y$ lies in $\RR^2$ seen as the complex plane, there is no control, and the dynamic is given by $f(y)= i\ y$, where $i^2=-1$. We clearly have
$$V_{\theta^k}(y_0)\xrightarrow[k\to\infty]{} \frac{1}{2\pi}\int_0^{2\pi} g(|y_0|e^{rit}) \deriv t,$$
for any sequence of evaluations $(\theta^k)_k$ satisfying the LTC.
\end{exa}

\begin{exa}\label{example:2} $y$ lies in the complex plane again, with $f(y,u)=i\ y\ u$, where $u\in U$ is a given bounded subset of $\RR$, and $g$ is any continuous function in $y$ (which thus does not depend on $u$).
\end{exa}

\begin{exa}\label{example:3} $f(y,u)=-y+u$, where $u\in U$ a given bounded subset of $\RR^d$, and $g$ is any continuous function in $y$ (which thus does not depend on $u$).
\end{exa}

\noindent We will show later (using Corollary \ref{cor:nonexpan}) that the general limit value exists in Examples \ref{example:2} and \ref{example:3}.

\section{On the long-term condition (LTC)} \label{sec:LTC}
In this section, we discuss the LTC. First, we give the following remarks.

\begin{rem} \label{rem:LTC}  (a). By definition, one has  $$\label{eq1} \forall s \geq 0, \forall t \geq 0, \forall \theta \in \De(\RR_+), \ TV_{s+t}(\theta)\leq TV_s(\theta)+TV_t(\theta).$$ This implies that $(\theta^k)_{k\geq 1}$ satisfies the LTC if and only if $ \overline{TV}_{1}(\theta^k)\xrightarrow[k\to\infty]{} 0$.\\
(b). If one takes $Q=\RR_+$ in definition of $TV_s(\theta^k)$ for each $s\geq 0$ and each $k\geq 1$,
we deduce that if $(\theta^k)_{k\geq 1}$ satisfies the LTC, then $\theta^k ([0,s]) \xrightarrow[k\to\infty]{}0$ for any  $s \geq 0$.
\end{rem}

\begin{rem}  \label{rem:LTC_a.c.} Let $\theta$ be an evaluation absolutely continuous $w.r.t.$ the Lebesgue measure on $\RR_+$, and $f_\theta$ its density.  Scheff\'e Theorem (cf. \cite {Devroye_85}, Theorem 1 in p.2) implies that: $$\forall s\geq 0,\ 2TV_s(\theta)=I_s(\theta)=_{def}\int_{0}^\infty|f_{\theta}(t+s)-f_{\theta}(t)|\deriv t.$$ Thus, if $(\theta^k)_{k\geq 1}$ is a sequence of evaluations with densities $(f_{\theta^k})_{k\geq 1}$:\\
(a). $(\theta^k)_{k\geq1}$ satisfies the LTC if and only if $\sup_{0\leq s\leq 1}I_s(\theta^k)\xrightarrow[k\to\infty]{}0$. If moreover, for each $k\geq 1$, $t\mapsto f_{\theta^k}(t)$ is non increasing on $\RR_+$, then $(\theta^k)_{k\geq 1}$ satisfies the LTC if and only if $\forall s\geq 0$, $\theta^k([0,s])=\int_{t=0}^\infty f_{\theta}(t)\deriv t - \int_{t=0}^\infty f_{\theta}(t+s)\deriv t \xrightarrow[k\to\infty]{} 0$.\\
(b). if $(\theta^k)_{k\geq1}$ satisfies the LTC, then $\int_{t=0}^{\infty}tf_{\theta^k}(t)\deriv t\xrightarrow[k\to\infty]{}\infty$. Indeed, Chebychev's inequality gives that $\int_{t=0}^{\infty} t f_{\theta^k}(t)\deriv t\geq M\left(1-\theta^k([0,M])\right)$ for all $M>0$.
\end{rem}

Here we discuss several cases where the LTC condition is satisfied.
\begin{exa} (Uniform distributions) Assume that for each $k$, $\theta^k$ is the uniform law over the interval $[a_k, b_k]$, with $0\leq a_k \leq b_k$. For each $k$,
\begin{itemize}
 \item \underline{$s\geq b_k-a_k$}: \  $I_s(\theta^k)=\left\{\begin{array} {ccc}
\frac{2}{b_k-a_k} & \mbox{if} & 0<s<a_k \\
\frac{1+(b_k-s)}{b_k-a_k} & \mbox{if} &  a_k<s<b_k \\
\frac{1}{b_k-a_k} & \mbox{if} &  b_k<s \\
\end{array}\right.$,

\item \underline{$s<b_k-a_k$}: \  $I_s(\theta^k)=\left\{\begin{array} {ccc}
\frac{2s}{b_k-a_k} & \mbox{if} & 0<s<a_k \\
\frac{2s}{b_k-a_k} & \mbox{if} &  a_k<s<b_k \\
\end{array}\right.$.
\end{itemize}

\noindent One can check easily that $(\theta^k)_k$ satisfies the LTC if and only if $b_k-a_k \xrightarrow[k\to\infty]{} \infty$. Indeed, by Remark \ref{rem:LTC_a.c.} (a), it is sufficient to look at $I_s(\theta^k)$ for $s\in[0,1]$.
\end{exa}

\begin{exa}\label{exaexpon}
(Abel average) Assume that for each $k$, $\theta^k$ has density $s\mapsto f_{\theta^k}(s)=\lambda_k e^{-\lambda_k s}\mathds{1}_{\RR_+}(s)$, with $\lambda_k>0$. Since $\forall k\geq 1$, $s\mapsto f_{\theta^k}(s)$ is non increasing, Remark \ref{rem:LTC_a.c.} (a) implies that $(\theta^k)_k$ satisfies the LTC if and only if:
$\forall T>0, \ \theta^k([0,T])=\int_{s=0}^T\lambda_k e^{-\lambda_k s}\deriv s=1-e^{-T/{\lambda_k}}\xrightarrow[k\to\infty]{} 0,$ which is again equivalent to $\lambda_k \xrightarrow[k\to\infty]{} 0.$
\end{exa}

\begin{exa}\label{exanormal}
(Folded normal distributions) Assume that for each $k$, $\theta^k$ is the distribution of a random variable $|X^k|$, where $X^k$ follows a normal law ${\cal N}(m_k, \sigma_k^2)$. The density of $\theta^k$ is given by:
$$\forall t \geq0, \; f_{\theta^k}(t) = \frac{1}{\sigma_k \sqrt{2\pi}} \left[\exp\left(-\frac{1}{2} {\left(\frac{t-m_k}{\sigma_k}\right)}^2 \right)+ \exp\left(-\frac{1}{2} {\left(\frac{t+m_k}{\sigma_k}\right)}^2\right)\right].$$
\end{exa}
\begin{claim}\label{claim:folded_normal}
$(\theta^k)_k$ satisfies the LTC if and only if $\sigma_k \xrightarrow[k\to\infty]{}\infty$.
\end{claim}
Our argument relies on the following lemma, whose proof is put in the \textbf{Appendix}.
Without loss of generality, we may assume that $m_k$ is non-negative for each $k$.
\begin{lemme}\label{lem:variation}
Let $\theta$ be the distribution of $X$ where $|X|$ follows the normal law $\N(m,\sigma)$ with $m,\sigma>0$. There exists some $t^*\in[0,m)$ such that $f'_\theta(t)>0$ for any $t\in(0,t^*)$ and  $f'_\theta(t)<0$ for any $t\in(t^*,\infty)$. Moreover, such $t^*$ satisfies that: $(t^*)^2\geq m^2-\sigma^2$.
\end{lemme}
\noindent \textbf{Proof of Claim \ref{claim:folded_normal}} \ We apply Lemma \ref{lem:variation} to each evaluation $\theta_k$ to obtain some $t^*_k\in[0, m_k)$ such that: $f_{\theta^k}(\cdot)$ is increasing on $[0,t^*_k)$ and decreasing on $[t^*_k, \infty)$. This enables us to write: $$\forall s\leq t^*_k, \ \ I_s(\theta^k)=\int_{t^*_k-s}^{t^*_k}f_{\theta^k}(t)\deriv t+\int_{t_k^*-s}^{t^*_k}|f_{\theta^k}(t+s)-f_{\theta^k}(t)|\deriv t+\int_{t^*_k}^{t^*_k+s}f_{\theta^k}(t)\deriv t.$$
We deduce then $s f_{\theta^k}(t^*_k-s)\leq I_{\theta^k}(s)\leq 4s f_{\theta^k}(t^k_*)$ for $s\leq t_k^*$. Assume below $\hat{t}^*=_{def}\liminf_{k\to\infty} t^*_k>0$, and the analysis is analogue for $\hat{t}^*=0$, which we omit here.\\
\noindent (*). Suppose that $\sigma_k\to\infty$, then $$f_{\theta^k}(t^*_k)=\frac{1}{\sigma_k \sqrt{2\pi}} \left[\exp\left(-\frac{1}{2} {\left(\frac{t^*_k-m_k}{\sigma_k}\right)}^2\right)+ \exp\left(-\frac{1}{2} {\left(\frac{t^*_k+m_k}{\sigma_k}\right)}^2\right)\right]\leq\frac{2}{\sigma_k\sqrt{2\pi}}\xrightarrow[k\to\infty]{}0.$$
This implies that for $S=\hat{t}^*\wedge 1$, \  $\sup_{0\leq s\leq S}I_s(\theta^k)\xrightarrow[k\to \infty]{} 0$.\\
\noindent (**). Conversely, suppose that $(\theta^k)_k$ satisfies the LTC. Then for any $s<\hat{t}^*$, $I_s(\theta^k)$ thus $f_{\theta^k}(t^*_k-s)$ vanishes as $k$ tends to infinity. This implies that either $\sigma_k\to \infty$ or $(\sigma_k)_k$ is bounded and $\big(m_k-(t^*_k-s)\big)_k\to \infty$. Lemma \ref{lem:variation} shows that the specified point $t^*_k$ for the evaluation $\theta_k$ satisfies $(t^*_k)^2\geq m_k^2-\sigma^2_k$, thus $m_k-(t^*_k+s)\leq m_k-t^*_k\leq \frac{\sigma^2_k}{m_k+t^*_k}\leq \frac{\sigma^2_k}{m_k}$. If $(\sigma_k)_k$ is bounded, $(m_k-t^*_k)_k$ thus $(m_k)_k$ should tend to infinity, but this leads to a contradiction with $m_k-t^*_k\leq \frac{\sigma^2_k}{m_k}$. $\hfill \Box$

Now we link the LTC condition to the discrete time framework. In a discrete time dynamic optimization problem, a general evaluation on the payoff stream is a probability distribution over $\NN^*=\NN/\{0\}$ the set of postive integers. For any $\xi=(\xi_1,...,\xi_t,...)$ in $\De(\NN^*)$, its "total variation" $TV(\xi)=\sum_{m=1}^\infty|\xi_{m+1}-\xi_m|$ is the stage by stage absolute difference between the measure $\xi$ and its one-stage "shift" measure $\xi'=(\xi_2,...,\xi_{t+1},...)$. (cf. Sorin \cite{Sorin_2002} or Renault \cite{Renault_2013}).

When the sequence of evaluations in continuous time admits step functions as densities, this link to discrete time framework is much clearer as seen by the following
\begin{pro} Let $(\theta^k)_k$ be a sequence of absolutely continuous evaluations in $\De(\RR_+)$, and their densities are given as: $\forall k\geq 1, f_{\theta^k} = \sum_{m=1}^\infty \xi^k_m {\mathds{1}}_{[m-1,m)}$, where $\xi^k=(\xi^k_1,...,\xi^k_m,...,)\in\De(\NN^*)$. Then $(\theta^k)_{k}$ satisfies the LTC if and only if $\sum_{m=1}^\infty | \xi^k_{m+1}-\xi^k_m|\xrightarrow[k\to\infty]{}0.$
\end{pro}
\noindent \textbf{Proof}: \ Fix $s\in[0,1]$. We shall write for each $k$,  $$I_s(\theta^k)=\sum_{m=1}^\infty\int_{[m-1,m)} \Big|f_{\theta^k}(t+s)-f_{\theta^k}(t)\Big|\deriv t.$$ For each $m=1,2,...$, we have
\begin{align*}
\int_{[m-1,m)} \Big|f_{\theta^k}(t+s)-f_{\theta^k}(t)\Big|\deriv t=&\int_{[m-1,m-s)}\Big|f_{\theta^k}(t+s)-f_{\theta^k}(t)\Big|\deriv t+\int_{[m-s,m)}\Big(f_{\theta^k}(t+s)-f_{\theta^k}(t)\Big)\deriv t\\
=& s \big|\xi^k_{m+1}-\xi^k_m\big|.
\end{align*}
As a consequence, $$I_s(\theta^k)=s\sum_{m=1}^\infty\big|\xi^k_{m+1}-\xi^k_m\big|\leq \sum_{m=1}^\infty\big|\xi^k_{m+1}-\xi^k_m\big|, \ \forall s\in[0,1].$$ In view of Remark \ref{rem:LTC_a.c.}, $(\theta^k)_k$ satisfies the LTC if and only if $\sum_{m=1}^\infty | \xi^k_{m+1}-\xi^k_m|\xrightarrow[k\to\infty]{} 0.$ $\hfill\Box$

We end this section by a preliminary lemma, which will be useful in later results.
\begin{lem} \label{lem:hahn}
Fix any $\theta\in\De(\RR_+)$ and any $t\in\RR_+$, we have
$$\left|\int_{[0,+\infty)}h(s)\deriv\theta(s)- \int_{[0,+\infty)} h(s-t)\deriv\theta(s)\right|\leq TV_t(\theta)\text{ and } \left|\int_{[0,+\infty)} h(s)\deriv\theta(s)- \int_{[0,+\infty)} h(s+t)\deriv\theta(s)\right|\leq 2TV_t(\theta),$$

\noindent for any $h(\cdot)\in\mathcal{M}(\RR_+,[0,1])$, where $\mathcal{M}(\RR_+,[0,1])=\Big\{h(\cdot)\Big|h:\RR_+\to[0,1], \text{Borel measurable}\Big\}$.
\end{lem}
\noindent \textbf{Proof}: \ We fix any $\theta\in\De(\RR_+)$ and $t\in\RR_+$. By definition of $\T_s\sharp\theta$, we have that for any $h(\cdot) \in\mathcal{M}(\RR_+,[0,1])$:
\begin{eqnarray}\label{eq:hahn1} \int_{[0,+\infty)} h(s)\deriv\theta(s)- \int_{[t,+\infty)} h(s-t)\deriv\theta(s)=\int_{[0,+\infty)} h(s)\deriv\theta(s)- \int_{[0,+\infty)} h(s)\deriv\T_{-t}\sharp\theta(s)
\end{eqnarray}
 and
\begin{eqnarray} \label{eq:hahn2} \int_{[0,+\infty)} h(s)\deriv\theta(s)- \int_{[0,+\infty)} h(s+t)\deriv\theta(s)=\int_{[0,+\infty)} h(s)\deriv\theta(s)- \int_{[0,+\infty)} h(s)\deriv\T_t\sharp\theta(s).
\end{eqnarray}
Since $\T_{-t}\sharp\theta$ and $\T_{t}\sharp\theta$ are both Borel measures on $\RR_+$, "$\theta-\T_{-t}\sharp\theta$" and "$\theta-\T_{t}\sharp\theta$" are both signed measures. Hahn's decomposition theorem\footnote{The first author acknowledges Eilon Solan for the discussion on using Hahn's decomposition theorem.} implies that:
\begin{eqnarray*}  \sup_{h\in \mathcal{M}(\RR_+,[0,1])}\left|\int_{[0,+\infty)} h(s)\deriv\theta(s)-\int_{[0,+\infty)} h(s)\deriv \T_{-t}\sharp \theta(s)\right|=\sup_{Q\in\mathcal{B}(\RR_+)}\Big|\theta(Q)-\T_{-t}\sharp \theta(Q)\Big|.
\end{eqnarray*}
and
\begin{eqnarray*} \sup_{h\in \mathcal{M}(\RR_+,[0,1])}\left|\int_{[0,+\infty)} h(s)\deriv\theta(s)-\int_{[0,+\infty)} h(s)\deriv \T_{t}\sharp \theta(s)\right|=\sup_{Q\in\mathcal{B}(\RR_+)}\Big|\theta(Q)-\T_{t}\sharp \theta(Q)\Big|.
\end{eqnarray*}
Combining with (\ref{eq:hahn1})-(\ref{eq:hahn2}), we obtain:
$$\left|\int_{[0,+\infty)} h(s)\deriv\theta(s)-\int_{[t,+\infty)} h(s-t)\deriv\theta(s)\right|\leq \sup_{Q\in\mathcal{B}(\RR_+)}\Big|\theta(Q)-\theta(Q+t)\Big|=TV_t(\theta)$$
and
\begin{align*}\left|\int_{[0,+\infty)} h(s)\deriv\theta(s)-\int_{[0,+\infty)} h(s+t)\deriv\theta(s)\right|\leq \sup_{Q\in\mathcal{B}(\RR_+)}\Big|\theta(Q)-\theta(Q-t)\Big|\leq \theta\big([0,t)\big)+ TV_t(\theta)\leq 2TV_t(\theta).
\end{align*}
The proof of the lemma is complete. $\hfill\Box$

\section{Main Result}\label{sec:result}

As will be shown in our main result, the function $V^*(y_0)$ defined in (\ref{eq:V*}) characterizes the general limit value of the optimal control problem  in case of convergence. We first rewrite it as
\begin{align*}
V^*(y_0)=\sup_{\mu \in \De(\RR_+)}\inf_{t \in \RR_+} V_{\T_t \,\sharp\,  \mu}(y_0)=\sup_{\mu\in \De(\RR_+)}\inf_{\overline{y}\in R(y_0)}V_{\mu}(\overline{y}).
\end{align*}
We give the following interpretation: consider the auxiliary optimal control problem (game) where an adversary of the controller chooses an evaluation $\mu$, and then knowing $\mu$ as given, the controller chooses an initial state in the reachable set $R(y_0)$. The running cost from time $t$ is evaluated by $\mu$ and $V^*(y_0)$ is the value of this problem starting from $y_0$.

\noindent Recall that a metric space $X$ is \textit{totally bounded} if for each $\varepsilon>0$, $X$ can be covered by finitely many balls of radius $\varepsilon$.
\begin{thm}\label{thm:main} Let $(\theta^k)_{k\geq 1}$ be a sequence of evaluations satisfying the LTC. Assume (\ref{ex:assumptions}) for the optimal control problem $\mathcal{J}=\langle U,g,f \rangle$. Then,\\
(i). $V^*=\sup_{k\in \NN}\inf_{t\in \RR_+} V_{\T_t \, \sharp \, \theta^k}.$\\
(ii). Any accumulation point (for the uniform convergence) of the sequence $(V_{\theta^k})_{k}$ is equal to $V^*$.\\
(iii). The sequence $(V_{\theta^k})_{k}$ uniformly converges if and only if the space $(\{V_{\theta^k}\},||\cdot||_\infty)$ is totally bounded.
\end{thm}

\begin{rem} \label{rem:sub-LTC} Let $(\theta^k)_k$ be a sequence of evaluations which contains a subsequence $(\theta^{\varphi_k})_k$ satisfying the LTC. Then Part $(i)$ of Theorem \ref{thm:main} still holds true for $(\theta^k)_k$.
\end{rem}

A more precise convergence result is obtained if we suppose that there exists a compact set $Y\subseteq\RR^d$ which is \textit{invariant} for the dynamic (\ref{eq:dynamics}), $i.e.$, such that $y(t,u,y_0)\in Y$ for all $u(\cdot)\in\mathcal{U}$, $t\geq 0$ and $y_0$ in $Y$.

\begin{defi}\label{defi:uniform-convergence}  Let $V$ be a function defined on $\RR^d$. In the optimal control problem $\mathcal{J}$, there is the \textbf{general uniform convergence} of the value functions $\{V_\theta\}$ to $V$ if:
$$\forall \varepsilon>0, \exists S>0, \exists \eta>0 \  s.t. \ \forall \theta\in\De(\RR_+), \text{with} \ \overline{TV}_S(\theta)\leq \eta, \ ||V_\theta-V||_\infty\leq \varepsilon. $$
\end{defi}

\begin{lem}\label{lem:uniform-convergence} Let $V$ be a function defined on $\RR^d$. The optimal control problem $\mathcal{J}$ admits $V$ as the general limit value if and only if there is the general uniform convergence of the value functions $\{V_\theta\}$ to $V$.
\end{lem}
\begin{proof} The general uniform convergence of the value functions $\{V_\theta\}$ to $V$ implies the existence of general limit value given as $V$.  Next we show that the existence of the general limit value given as $V$ is sufficient to deduce the general uniform convergence of $\{V_\theta\}$ to $V$. Suppose by contradiction that there is no general uniform convergence of $\{V_\theta\}$ to $V$, $i.e.$,
\begin{align*}
\exists \varepsilon_0>0, \ \forall S>0,  \ \forall \eta^k>0, \ \exists \theta^k \in\De(\RR_+) \text{ with } \overline{TV}_S(\theta^k)\leq \eta^k, \text{ and } ||V_{\theta^k}-V||_\infty>\varepsilon_0, \ \forall k\geq 1.
\end{align*}
Let $\varepsilon_0>0$ be fixed as above. We take a vanishing positive sequence $(\eta^k)_k$ and some $S_0>0$, then there is a sequence of evaluations $(\theta^k)$ with $\overline{TV}_{S_0}(\theta^k)\leq \eta^k\xrightarrow[k\to\infty]{}0$, and $\liminf_{k}||V_{\theta^k}-V||_\infty\geq \varepsilon_0$.  According to Remark \ref{rem:LTC} (a), such $(\theta^k)_k$ satisfies the LTC, while $(V_{\theta^k})$ does not converges uniformly to $V^*$. This is a contradiction. \end{proof}

\begin{cor}\label{cor:uniformcontin} Assume (\ref{ex:assumptions}) for the optimal control problem $\mathcal{J}=\langle U,g,f \rangle$. Suppose that there is a compact set $Y\subseteq \RR^d$ which is invariant for the dynamic (\ref{eq:dynamics}), and that the family $\{V_{\theta}:\ \theta\in\De(\RR_+)\}$ is uniformly equicontinuous on $Y$. Then there is the general uniform convergence of the value functions $\{V_\theta\}$ to $V^*$.
\end{cor}

\noindent \textbf{Proof}: \  By assumption, the family of value functions $\{V_\theta: \theta\in\De(\RR_+)\}$ is both uniformly bounded and uniformly equicontinuous on the compact invariant set $Y$, so we can use Ascoli's theorem to deduce the totally boundedness of the space $(\{V_\theta\},||\cdot||_\infty)$. Theorem \ref{thm:main} implies that: for any $(\theta^k)_k$ satisfying the LTC, the corresponding sequence of value functions $(V_{\theta^k})$ converges uniformly to $V^*$ as $k$ tends to infinity. Thus $\mathcal{J}$ has a general limit value given as $V^*$, and according to Lemma \ref{lem:uniform-convergence}, there is the uniform convergence of value functions $\{V_\theta\}$ to $V^*$. $\hfill\Box$

\noindent We shall give the existence result of the general limit value under sufficient conditions expressed directly in terms of properties of the control dynamic (\ref{eq:dynamics}) and of the running cost $g$.

Let us introduce the following \textit{non expansive} condition (cf. \cite{Quincampoix_2011}). The control dynamic (\ref{eq:dynamics}) is non expansive if $$\forall y_1,y_2\in\RR^d,\ \sup_{a\in U}\inf_{b\in U}\Big\langle y_1-y_2,f(y_1,a)-f(y_2,b)\Big\rangle\leq 0.$$

\begin{defi}\label{def:compact-nonexpan} The optimal control problem $\mathcal{J}=\langle U,g,f \rangle$ is called \textbf{compact non expansive} if it satisfies the following three conditions: \\
$(A.1)$ there is a compact set $Y\subseteq \RR^d$ is the invariant for the dynamic (\ref{eq:dynamics});\\
$(A.2)$ the running cost function $g(\cdot)$ does not depend on $u\in U$, and is continuous in $y\in\RR^d$;\\
$(A.3)$ the control dynamic (\ref{eq:dynamics}) is non expansive on $Y$.
\end{defi}

\begin{cor}\label{cor:nonexpan} Assume (\ref{ex:assumptions}) for the optimal control problem $\mathcal{J}=\langle U,g,f \rangle$. Suppose that that $\mathcal{J}$ is compact non expansive, then the general limit value exists in $\mathcal{J}$ and is given as $V^*$.
\end{cor}
\noindent \textbf{Proof}: \ Under $(A.1)$ and $(A.3)$, Proposition 3.7 in \cite{Quincampoix_2011} implies that:
\begin{align}\label{eq:nonexpansive}
 \forall (y_1,y_2) \in Y^2, \forall u(\cdot)\in\mathcal{U}, \exists v(\cdot)\in\mathcal{U},
s.t. \ \forall t\geq 0, \  ||y(t,u,y_1)-y(t,v,y_2)||\leq ||y_1-y_2||.
\end{align}
\noindent We claim that the family $(V_{\theta})_{\theta\in\De(\RR_+)}$ is uniformly equicontinuous on $Y$, thus Corollary \ref{cor:uniformcontin} and Lemma \ref{lem:uniform-convergence} apply. Fix any $(y_1,y_2)\in Y^2$, $\theta\in\De(\RR_+)$, and $\varepsilon>0$. Let $u$ be $\varepsilon$-optimal for $V_\theta(y_1)$: $$V_\theta(y_1)\geq \int_{[0,+\infty)} g\Big(y(s,u,y_1)\Big)\deriv\theta(s)-\varepsilon.$$
According to the non expansive property, there exists $v(\cdot)$ in $\mathcal{U}$ as in (\ref{eq:nonexpansive}) such that
\begin{align}\label{eq:NC}
||y(s,u,y_1)-y(s,v,y_2)||\leq ||y_1-y_2||, \ \forall s\geq 0.
\end{align}
By definition, $V_\theta(y_2)\leq \int_{[0,+\infty)}g\big(y(s, y_2,v)\big)\deriv\theta(s)$,
hence \begin{eqnarray*} V_{\theta}(y_2)-V_{\theta}(y_1)\leq \int_{[0,+\infty)} \Big[g\big(y(s,v,y_2))-g(y(s,u,y_1)\big)\Big]\deriv\theta(s)+\varepsilon.\end{eqnarray*}
\noindent Denoting $\omega_g$ the modulus of continuity of $g$, we obtain in view of ($\ref{eq:NC}$):
\begin{eqnarray*}
V_{\theta}(y_2)-V_{\theta}(y_1)&\leq& \int_{[0,+\infty)} \Big[g\big(y(s,v,y_2)\big)-g\big(y(s,u,y_1)\big)\Big]\deriv\theta(s)+\varepsilon\leq \omega_g(||y_1-y_2||)+\varepsilon.
\end{eqnarray*}
Interchanging $y_1$ and $y_2$ and taking into account of $\varepsilon>0$ being arbitrary, we deduce that  $(V_{\theta})_{\theta\in\De(\RR_+)}$ is uniformly equicontinuous on the invariant set $Y$. This finishes the proof. $\hfill\Box$

\begin{rem} Both Example \ref{example:2} and Example \ref{example:3} satisfy conditions of Corollary \ref{cor:nonexpan}, so there is general uniform convergence of the value functions $\{V_\theta\}$ (the existence of the general limit value).
\end{rem}
\begin{rem} Our result generalizes Proposition 3.3 in \cite{Quincampoix_2011} which proved the uniform convergence of the $t$-horizon values in compact non expansive optimal control problems.
\end{rem}

We end this section by presenting two (counter)examples, showing that the results in Theorem \ref{thm:main} do not hold if some of their conditions is not satisfied.

\noindent The first example is an uncontrolled dynamic. We show that if $(\theta^k)_k$ contains no subsequence satisfying the LTC, then the result in Part (i) of Theorem \ref{thm:main} does not hold, $i.e.$, $\sup_{k\geq 1}\inf_{t\geq 0} V_{\T_t\sharp \theta^k} (y_0)<\sup_{\theta\in\De(\RR_+)}\inf_{t\geq 0} V_{\T_t\sharp \theta} (y_0)$ for some $y_0$ (cf. Remark \ref{rem:sub-LTC}).

\begin{countexa}\label{countexa} Consider the uncontrolled dynamic on $\RR$: $y(0)=y_0$ and $y'(t)=-\left(y(t)-1\right), \forall t\geq 0$. The trajectory is then $y(t)=1+(y_0-1)e^{-t}$. The running cost function $g:\RR\to [0,1]$ is given by:
$$g(y)=\left\{ \begin{array} {ccc}
 0& \mbox{if} &    y < 0 \\
 y &\mbox{if} &  0\leq y\leq 1 \\
1 & \mbox{if } & y > 1\\\end{array} \right.$$
\noindent We have that $V^*(y_0)=\sup_{\theta \in \De(\RR_+)}\inf_{t \in \RR_+} V_{\T_t\sharp\theta}(y_0)=1, \ \forall y_0\in\RR$. Indeed, let $y_0$ be given and fix any $\varepsilon>0$, there is some $T_\varepsilon>0$ such that $|y(T)-1|\leq \varepsilon$ for all $T\geq T_\varepsilon$. Take an evaluation $\theta$ in $\De(\RR_+)$ with $\theta([0,T_\varepsilon])=0$. This enables us to deduce that: for all $t\geq 0$,  $$V_{\T_t\sharp\theta}(y_0)=\int_{[T_\varepsilon,+\infty)} g\big(y(s+t)\big)\deriv\theta(s)\geq \int_{[T_\varepsilon, +\infty)} g\big(y(T_\varepsilon)\big)\deriv\theta(s)\geq \big(1-y(T_\varepsilon)\big)\theta([T_\varepsilon,\infty])\geq 1-\varepsilon.$$

\begin{tikzpicture}
\centering
\draw [thick, <->] (6,0) -- (0,0) -- (0,3);
\node [below] at (6.2,0) {\text{time t}};
\node [above] at (0,3) {\text{distance of $y(t)$ from 1}};
\node [below] at (-0.8,2.6) {$|y_0-1|$};
\draw [thin] (-0.1,2.3)--(0.1,2.3);
\node [below] at (0,0) {\text{0}};
\draw [ultra thick, domain=0:5.5] plot (\x, {2.3*exp(-\x)+0.05});
\node [below] at (-0.3, 0.3) {$\varepsilon$};
\draw [dashed] (0,0.1) -- (3.5,0.1);
\draw [dashed] (3.5,1.0) -- (3.5,-0.10);
\node [below] at (3.5,0) {$T_\varepsilon$};
\end{tikzpicture}

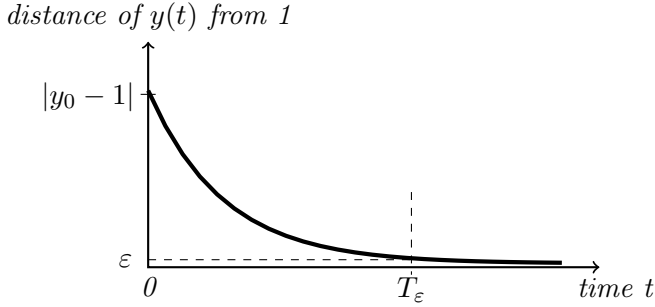
\noindent \textbf{Figure 5.1}: The solution $y(t)=1+(y_0-1)e^{-t}$ to the dynamic is represented by the thick curve. For given $\varepsilon>0$, $T_\varepsilon>0$ is chosen such that $|y(T_\varepsilon)-1|=\varepsilon$.\\

\noindent Consider now any sequence of evaluations $(\theta^k)_k$  which does not contain any subsequence satisfying the LTC. Under the assumption that the density $f_{\theta^k}$ for each evaluation $\theta^k$ is non increasing, we show that Part (i) of Theorem \ref{thm:main} is not valid: $V^*\neq \sup_{k\in\NN}\inf_{t\in\RR_+} V_{\T_t\sharp\theta^k }$.

\noindent Indeed, let us take any $y_0<1$ and suppose that $\sup_{k\in\NN}\inf_{t\in\RR_+} V_{\T_t\sharp \theta^k}(y_0)=V^*(y_0)$, which is equal to $1$ as was proved. Let $\varphi(k)$ be a subsequence such that $\lim_{k\to\infty} \inf_{t\in \RR_+} V_{\T_t\sharp \theta^{\varphi(k)}}(y_0)=1$. $(\theta^{\varphi(k)})_k$ does not satisfy the LTC by assumption, so Remark \ref{rem:LTC_a.c.} (a) implies that there exists some $T>0$ with $\theta^{\varphi(k)}([0,T])\nrightarrow0$. Let $\varphi_m$ be the subsequence of $\varphi$ and $\eta>0$ such that $\theta^{\varphi_m(k)}([0,T])\xrightarrow[k\to\infty]{}\eta$. We obtain for any $k\geq 1$, \begin{eqnarray*}\inf_{t\in\RR_+}V_{\T_t\sharp\theta^{\varphi_m(k)}}(y_0)\leq V_{\theta^{\varphi_m(k)}}(y_0)&=&\int_{[0,T]}g\left(y(t)\right)\deriv\theta^{\varphi_m(k)}(t)+\int_{[T,+\infty]}g\left(y(t)\right)\deriv\theta^{\varphi_m(k)}(t)\\
&\leq& y(T)\theta^{\varphi_k(m)}([0,T])+\theta^{\varphi_k(m)}([T,\infty]).
\end{eqnarray*}
This implies that for such fixed $y_0<1$ , $\lim_k\inf_{t\in\RR_+}V_{\T_t\sharp \theta^{\varphi_m(k)}}(y_0)\leq y(T)\eta+(1-\eta)<1$. This contradicts the assumption that $\sup_{k\in\NN}\inf_{t\in\RR_+}V_{\T_t\sharp\theta^k}(y_0)=1$, and our claim is proved.
\end{countexa}

In the second example, we study the convergence of the value functions of a control problem along two different sequences of evaluations satisfying the LTC. Along the first sequence, the value functions converge uniformly to $V^*$; while along the second, the value functions point-wisely converge, but not uniformly (thus the family of value functions is not totally bounded for the uniform norm), to a limit function which is different from $V^*$.

\begin{countexa}\label{exa:t-horizon}
Consider the control problem on the state space $\RR=(-\infty,+\infty)$, where the control set is $U=\{+1,-1\}$; the dynamic is\footnote {Notice that the dynamic is discontinuous at $y=0$ when $u=+1$. To get the desired asymptotic result under the Liptchitz regularity, one can slightly modify dynamic to set $f(y,+1)=y$ for $y\in[0,1]$ and others unchanged.}: $$f(y, u)=u \text{ for all }(y,u)\in\RR_+\times U \text{ and } f(y,u)=-1 \text{ for all } (y,u)\in\RR^*_{-}\times U,$$ where $\RR^*_{-}=\RR_-/\{0\}$; and the running cost function is:
$$g(y,u)=\left\{\begin{array}{ccc}
+1 & \mbox{if} & u=+1, \ y\geq 0\\
0 & \mbox{if} & u=-1, \ y\geq 0\\
+K & \mbox{if} & \ y<0     \\
\end{array}\right.$$

\noindent Suppose that $K>1$ large enough, so the cost on $\RR_-$ is positive and high. Whenever the state reaches $y=0$, it is optimal to choose control $u=+1$ and this drives the state back to $\RR_+$; on $\RR^*_{-}$, the dynamic is $f=-1$, independent of control and state. $V_\theta(y_0)=K$ for all $y_0$ in $\RR^*_{-}$ and $\theta$ in $\De(\RR_+)$, so the reduced state space is $\RR_+$, and we consider value functions defined on it.

\noindent $V^*(y_0)=\sup_{\theta}\inf_{t\geq 0} V_{\T_t\sharp\theta}(y_0)=0$ for any $y_0\geq 0$. Fix any $y_0\geq 0$. For any $\theta\in\De(\RR_+)$ and $\varepsilon>0$, let $t^\varepsilon\geq 0$ such that $\theta([0, t^\varepsilon])\geq 1-\varepsilon$. Define now the control $u^\varepsilon(\cdot)$ to be: $u^\varepsilon(t)=+1$, if $t\in[0, t^\varepsilon]$ and $u^\varepsilon(t)=-1$ if $t\in(t^\varepsilon, \infty)$, which gives: $\gamma_{\T_{t^\varepsilon}\sharp\theta}(y_0,u^\varepsilon)\leq \varepsilon K$.

\noindent Consider $(\theta^k)_k$ the sequence of evaluations with density $f_{\theta^k}(s)=\frac{1}{k}\mathds{1}_{[k,2k]}(s)$ for each $k$, and $(\bar{\theta}^k)_k$ the sequence of $k$-horizon evaluations with density $f_{\bar{\theta}^k}(s)=\frac{1}{k}\mathds{1}_{[0,k]}(s)$ for each $k$. We show that:

$(\{V_{\theta^k}\},||\cdot||_\infty)$ is totally bounded and $(V_{\theta^k})$ converges uniformly to $V^*$; while  $(\{V_{\overline{\theta}^k}\},||\cdot||_\infty)$ is not totally bounded and $(V_{\overline{\theta}^k})$ does not converge to $V^*$.

\noindent Let $y_0\geq 0$, we have that:
\begin{enumerate}
\item $V_{\theta^k}(y_0)=0$, for all $k\geq 1$. \ Indeed, one optimal control for $V_{\theta^k}(y_0)$ can be taken as: $u^*(t)=+1, \ t\in[0,k]$ and $u^*(t)=-1, \ t\in(k, 2k]$;
\item $V_{\overline{\theta}^k}(y_0)=0$ if $k\leq y_0$ and $V_{\overline{\theta}^k}(y_0)=\frac{1}{2}-\frac{y_0}{2k}$ if $k>y_0$. \ Indeed, for $k\leq y_0$, one optimal control for $V_{\overline{\theta}^k}(y_0)$ can be taken as: $u^*(t)=-1, \ t\in[0,k]$; for $k>y_0$, one optimal control for $V_{\overline{\theta}^k}(y_0)$ can be taken as: $u^*(t)=+1, \ t\in[0, \frac{k-y_0}{2}]$ and $u^*(t)=-1, \ t\in(\frac{k-y_0}{2},k]$, so $\gamma_{\overline{\theta}^k}(y_0,u^*)=\frac{(k-y_0)/2}{k}=\frac{1}{2}-\frac{y_0}{2k}$.
\end{enumerate}

See the following two pictures for illustration.

\begin{tikzpicture}

\centering
\node [below] at (0,-0.1) {\text{0}};
\node [below] at (8,-0.1) {\text{0}};
\draw [thick, <->] (6,0) -- (0,0) -- (0,3);
\draw [ultra thick] (0,1) --(2.5,3.5);
\draw [ultra thin] (2.5,3.5) -- (5,1);

\draw [dashed] (2.5,0) -- (2.5,3.5);
\draw [dashed] (0,1) -- (5,1);

\draw [decorate,decoration={brace,amplitude=7pt},xshift=0pt,yshift=0pt]
(0,0) -- (0,1) node [black,midway,xshift=-0.6cm]
{\footnotesize $y_0$};
\draw [decorate,decoration={brace,amplitude=10pt, mirror},xshift=0pt,yshift=0pt]
(2.5,0) -- (5,0) node [black,midway,yshift=-0.6cm]
{\footnotesize $k$};

\draw[thick] ([shift=(0:0.5)] 0,1) arc (0:45:0.5);
\node at (0.8, 1.3) {$\frac{\pi}{4}$};

\draw[thick] ([shift=(0:0.5)] 8,1) arc (0:45:0.5);
\node at (8.8, 1.3) {$\frac{\pi}{4}$};

\draw [thick, <->] (15.7,0) -- (8,0) -- (8,3);
\draw [decorate,decoration={brace,amplitude=7pt},xshift=0pt,yshift=0pt]
(8,0) -- (8,1) node [black,midway,xshift=-0.6cm]
{\footnotesize $y_0$};
\draw [ultra thick] (8,1) -- (11,4);
\draw [ultra thin] (11,4)--(15,0);
\draw [dashed] (8,1) -- (14,1);
\draw [dashed]  (11,4) -- (11,0);

\draw [decorate,decoration={brace,amplitude=10pt, mirror},xshift=0pt,yshift=0pt]
(8,0) -- (15,0) node [black,midway,yshift=-0.6cm]
{\footnotesize $k$};

\draw [decorate,decoration={brace,amplitude=10pt, mirror},xshift=0pt,yshift=0pt]
(11,1) -- (11,4) node [black,midway,yshift=-0.4cm]
{\footnotesize };
\node at (11.7,2.5) {$\frac{k-y_0}{2}$};

\node [below] at (2.5,-0.7) {\text{ one optimal control for } $\theta^k$};
\node [below] at (11.5,-0.7) {\text{ one optimal control for } $\overline{\theta}^k$, \ $k>y_0$};
\node [below] at (6.2,0) {\text{time}};
\node [below] at (15.7,0) {\text{time}};
\node [above] at (0,3) {\text{distance from 0}};
\node [above] at (8,3) {\text{distance from 0}};
\draw [ultra thick] (0,-2.2) -- (3,-2.2);
\node  [above] at (1.5,-2.2) {$u=+1$ \ and \ $g=1$};

\draw [ultra thin] (5,-2.2) -- (8,-2.2);
\node  [above] at (6.5,-2.2) {$u=-1$ \ and \ $g=0$};

\end{tikzpicture}

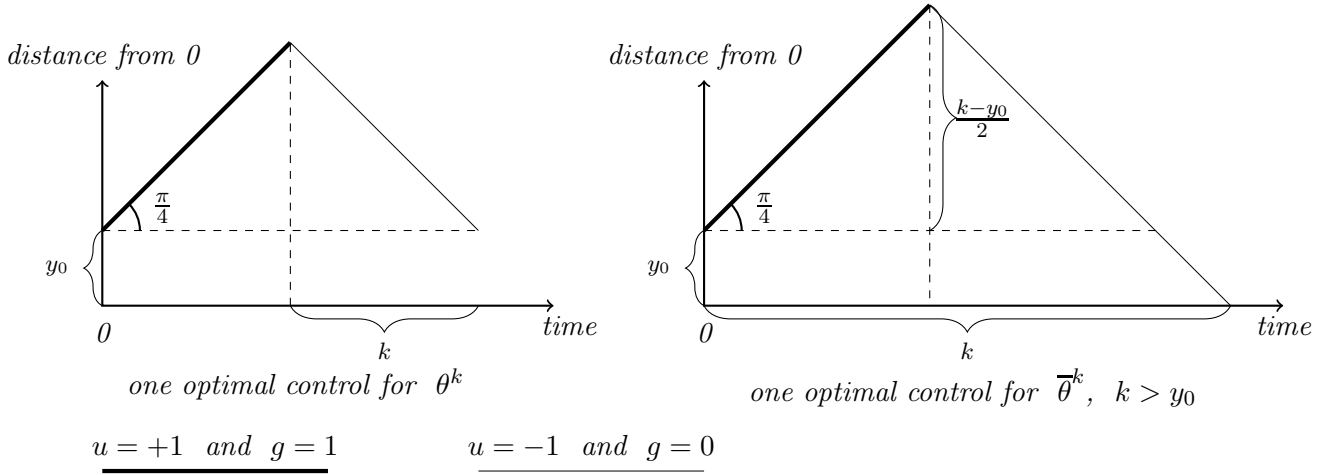
\noindent \textbf{Figure 5.2}: The left figure describes the dynamic of one optimal control for the evaluation $\theta^k$, which is $u^*=+1$ on $[0,k]$ and $u^*=-1$ on $(k,2k]$; the right figure describes the dynamic of one optimal control for the evaluation $\bar{\theta}^k$ with $k>y_0$, which is $u^*=+1$ on $[0,\frac{k-y_0}{2}]$ and $u^*=-1$ on $(\frac{k-y_0}{2},k]$. Here, the vertical axis represents the distance of $y(t)$ from zero and the thick trajectory ($resp.$ thin trajectory) corresponds to state on which $u=+1$ and $g=1$ ($resp.$ $u=-1$ and $g=0$).\\

\noindent We deduce that $\left(V_{\theta^k}(y_0)\right)_k$ converges uniformly to $V^*(y_0)=0$ on $\RR_+$; and that
$V_{\overline{\theta}^k}(y_0)\xrightarrow[k\to\infty]{} \frac{1}{2}$, while the convergence is not uniformly in $y_0\in\RR_+$: indeed, for all $k\geq 1$, $V_{\overline{\theta}^k}(k)=0$.
\end{countexa}

\section{Proof of main result: Theorem \ref{thm:main}}\label{sec:proof}

Consider in this section a sequence of evaluations $(\theta^k)_k$ that satisfies the LTC. As the proof is rather long, we divide it into two main parts:
\begin{itemize}
\item in Subsection \ref{sec:5.1}, we present the first preliminary result, Proposition \ref{propkey}. It is used in two ways: first, we obtain an immediate consequence of it for later use, which bounds $\liminf_k V_{\theta^k}$ from below in terms of the auxiliary value functions $\{V_{\T_t\sharp\theta^k}:k\in\NN^*, t\in\RR_+\}$; second, we deduce from it in Corollary \ref{coro:V*} the proof for Part (i) of Theorem \ref{thm:main}.
\item  In Subsection \ref{sec:5.2}, we prove Parts (ii)-(iii) of Theorem \ref{thm:main}. Lemma \ref{lem:limsup} gives an upper bound of $\limsup_k V_{\theta^k}$ in terms of the auxiliary value functions $\{V_{\T_t\sharp\theta^k}:k\in\NN^*, t\in\RR_+\}$, which is, together with the result from Proposition \ref{propkey}, used to end the proof.
\end{itemize}

\subsection{A first preliminary result and proof for Part (i)}\label{sec:5.1}

\begin{pro}\label{propkey}
For any $\mu$ in $\De(\RR_+)$, and any initial state $y_0$ in $\RR^d$,
$$\inf_{t\in\RR_+}V_{\T_t\sharp\theta}(y_0) \leq \liminf_{k} V_{\theta^k}(y_0).$$
In particular, we have for all $y_0$ in $\RR^d$,
$$\sup_{k\in\NN^*}\inf_{t\in\RR_+}V_{\T_t\sharp\theta^k}(y_0)\leq \liminf_k V_{\theta^k}(y_0).$$
\end{pro}

\noindent {\bf {Proof}:} \ \ Fixing $y_0$ and $\mu$, we set $\beta=_{def}\inf_{t\in\RR_+}V_{\T_t\sharp\theta}(y_0)$. For any $\varepsilon>0$ fixed, there exists some $T_0>0$ such that $\mu([T_0,\infty))< \varepsilon$.
Take any control $u(\cdot)$ in $\mathcal{U}$.
By definition of $\beta$, we have that $$\forall T\geq 0,\ \int_{[0,+\infty)} g\left(y(t+T,u,y_0),u(t+T)\right)\deriv\mu(t)\geq \beta,$$
\noindent thus
\begin{align}\label{eq:beta-1}
\forall T\geq 0,\ \int_{[0,T_0]} g\left(y(t+T,u,y_0),u(t+T)\right)\deriv\mu(t)\geq \beta-\varepsilon.
\end{align}
\noindent For each $k\geq 1$, integrating both sides of (\ref{eq:beta-1}) over $T\in[0,\infty)$ w.r.t. the evaluation $\theta^k$, we obtain
\begin{align}\label{eq:beta-2}
\int_{[0,+\infty)} \int_{[0,T_0]} g\left(y(t+T,u,y_0),u(t+T)\right)\deriv\mu(t)\deriv\theta^k(T)\geq \beta-\varepsilon.
\end{align}
Applying Fubini's Theorem to (\ref{eq:beta-2}) yields
\begin{align}\label{eq:beta-3}
\beta-\varepsilon\leq \int_{[0,T_0]}\left[\int_{[0,+\infty)} g\left(y(t+T,u,y_0),u(t+T)\right)\deriv\theta^k(T)\right]\deriv\mu(t)=\int_{[0,T_0]}\left[\gamma_{\T_t\sharp\theta^k}(y_0,u)\right]\deriv\mu(t),
\end{align}
where $\gamma_{\T_t\sharp \theta^k}(y_0,u)=\int_{[0,+\infty)} g\big(y(t+T,u,y_0),u(t+T)\big) \deriv\theta^k(T)$. According to Lemma \ref{lem:hahn}, we have $|\gamma_{\theta^k}(y_0,u)-\gamma_{\T_t\sharp\theta^k}(y_0,u)|\leq 2TV_{t}(\theta^k)$. This enables us to rewrite (\ref{eq:beta-3}) as:
\begin{eqnarray*}
\beta-\varepsilon &\leq& \int_{[0,T_0]} \left(\gamma_{\theta^k}(y_0,u)+2TV_t(\theta^k)\right)\deriv\mu(t)\\
&\leq& \left(\gamma_{\theta^k}(y_0,u)+2\overline{TV}_{T_0}(\theta^k)\right)\mu([0,T_0])\\
&\leq & \gamma_{\theta^k}(y_0,u)+2\overline{TV}_{T_0}(\theta^k).
\end{eqnarray*}
The control $u(\cdot)\in\mathcal{U}$ being taken arbitrarily, we deduce that
$$\beta-\varepsilon\leq V_{\theta^k}(y_0)+2\overline{TV}_{T_0}(\theta^k).$$
\noindent Since $(\theta^k)$ satisfies the LTC, $\overline{TV}_{T_0}(\theta^k)$ vanishes as $k$ tends to infinity. The proof is achieved. $\hfill \Box$

\noindent We end the proof for Part (i) of Theorem \ref{thm:main} by the following corollary of Proposition \ref{propkey}.

\begin{cor}\label{coro:V*} [Proof for Part (i) of Theorem \ref{thm:main}]
$$ \sup_{\mu\in\De(\RR_+)}\inf_{t\in\RR_+}V_{\T_t\sharp \mu}(y_0)=\sup_{k\geq 1}\inf_{t\in\RR_+} V_{\T_t\sharp\theta^k}(y_0),\ \forall y_0\in \RR^d.
$$
\end{cor}
\noindent{\bf{Proof}:} Fix $y_0\in\RR^d$, and denote $\varrho = \sup_{k\geq 1}\inf_{t\geq 0} V_{\T_t\sharp\theta^k} (y_0)$. It is clear that $\varrho \leq \sup_{\mu\in\De(\RR_+)}\inf_{t\geq 0} V_{\T_t\sharp\mu} (y_0)$. Now for the converse inequality, consider for each $k\geq 1$, there exists $m(k)$ in $\RR_+$ such that $V_{\T_{m(k)}\sharp\theta^k}(y_0)\leq \varrho+1/k$. Since $\T_{m(k)}\sharp \theta^k$ -- the image measure of $\theta^k$ by the function $s\mapsto s+m(k)$ -- is also an evaluation on $\RR_+$, we have:
$$\forall s\geq 0, \ TV_s(\T_{m(k)}\sharp\theta^k)=\sup_{Q\in\mathcal{B}(\RR_+)}\Big|\theta^k\big((Q-m(k))\cap\RR_+\big)-\theta^k\big((Q-m(k)+s)\cap\RR_+\big)\Big|\leq TV_s(\theta^k)+\theta^k([0,s]).$$
We deduce that $(\T_{m(k)}\sharp\theta^k)_k$ satisfies the LTC whenever $(\theta^k)_k$ does so. According to Proposition \ref{propkey},  $\forall \mu\in\De(\RR_+), \ \inf_{t\in\RR_+}V_{\T_t\sharp \mu}(y_0)\leq \liminf_k V_{\T_{m(k)}\sharp\theta^k}(y_0)\leq \varrho$, thus $\sup_{\mu\in\De(\RR_+)}\inf_{t\in\RR_+}V_{\T_t\sharp \mu}(y_0)\leq \varrho$. The proof is complete. $\hfill \Box$

\subsection{Proof for Parts (ii)-(iii)}\label{sec:5.2}
\noindent In this subsection, we give the proof for Parts (ii)-(iii) of Theorem \ref{thm:main}.
\noindent We begin with the following result, which compares the values under evaluation $\mu$ and its $t$-"shifted" evaluation $\T_t\sharp\mu$ for any $t>0$.
\begin{lem}\label{lemshift} Let $\mu$ in $\De(\RR_+)$ be any evaluation.
Then: for all $t\geq 0$ and $y_0\in \RR^d$,
$$ V_\mu(y_0)\leq V_{\T_t\sharp \mu}(y_0)+ 2TV_t(\mu).$$
\end{lem}
\noindent \textbf{Proof}: \ Fix $\mu\in\De(\RR_+)$, $t\geq 0$, $y_0\in\RR^d$. By Lemma \ref{lem:hahn}, we have
$$\gamma_\mu(y_0,u)\leq \gamma_{\T_t\sharp\mu}(y_0,u) + 2TV_t(\mu), \ \forall u(\cdot)\in\mathcal{U}.$$
\noindent For all $\varepsilon>0$, take $u^\varepsilon(\cdot)\in\mathcal{U}$ be an $\varepsilon$-optimal control for $V_{\T_t\sharp\mu}(y_0)$, $i.e.$, $\gamma_{\T_t\sharp\mu}(y_0,u^\varepsilon)\leq V_{\T_t\sharp\mu}(y_0)+\varepsilon$. We obtain that $$\gamma_\mu(y_0,u^\varepsilon)\leq V_{\T_t\sharp\mu}(y_0)+\varepsilon+ 2TV_t(\mu).$$
\noindent Since $V_\mu(y_0)=\inf_{u(\cdot)\in\mathcal{U}} \gamma_\mu(y_0,u)$ and $\varepsilon>0$ being arbitrary, we deduce that
$$ V_{\mu}(y_0)\leq V_{\T_t\sharp\mu}(y_0)+ 2TV_t(\mu),$$
which proves the lemma. $\hfill\Box$

\noindent The following result gives an upper bound on "$\limsup_kV_{\theta^k}$" in terms of the auxiliary value functions $\{V_{\T_t\sharp\theta^k}:k\in\NN^*, t\in\RR_+\}$.
\begin{lem}\label{lem:limsup} For all $T_0\geq 0$ and any $y_0$ in $\RR^d$,
$$ \limsup_k V_{\theta^k}(y_0)= \limsup_k \inf_{t\leq T_0} V_{\T_t\sharp\theta^k}(y_0).$$
In particular, for all $T_0\geq 0$ and any $y_0$ in $\RR^d$, $$\limsup_kV_{\theta^k}(y_0)\leq \sup_{k\in\NN^*}\inf_{t\leq T_0}V_{\T_t\sharp \theta^k}(y_0).$$ %
\end{lem}
\noindent{\bf Proof:} Fix $T_0\geq 0$ and $y_0\in\RR^d$. The inequality $''\limsup_{k}\inf_{t\leq T_0}V_{\T_t\sharp\theta^k}\leq \limsup_kV_{\theta^k}''$ is clear by taking $t=0$ for each $k$. Now for the converse inequality $''\limsup_{k}\inf_{t\leq T_0}V_{\T_t\sharp\theta^k}\geq \limsup_kV_{\theta^k}''$: according to Proposition \ref{lemshift}, we have that for all $k$ and $t\leq T_0$, $$V_{\theta^k}(y_0)\leq V_{\T_t\sharp \theta^k} (y_0)+ 2TV_{t}(\theta^k).$$
For each $k\geq 1$, take $t^k\leq T_0$ with $V_{\T_{t^k}\sharp \theta^k}(y_0)\leq \inf_{0\leq t\leq T_0}V_{\T_{t^k}\sharp\theta^k}+\frac{1}{k}$, which gives us:
\begin{eqnarray*}V_{\theta^k}(y_0)&\leq& \inf_{0\leq t\leq T_0}V_{\T_t\sharp \theta^k} (y_0)+\frac{1}{k}+2TV_{t^k}(\theta^k)\\ &\leq& \inf_{0\leq t\leq T_0}V_{\T_t\sharp \theta^k} (y_0)+\frac{1}{k}+ 2\overline{TV}_{T_0}(\theta^k).
\end{eqnarray*}
\noindent Since $(\theta^k)_k$ satisfies the LTC, $\overline{TV}_{T_0}(\theta^k)$ vanishes as $k$ tends to infinity. By taking "$\limsup_k$" on both sides of above inequality, the proof of the lemma is complete.
$\hfill \Box$

\noindent Now we end the proof for Theorem \ref{thm:main}. To do this, we first summarize results in Proposition \ref{propkey} and Lemma \ref{lem:limsup} in the following chain form, which is then used for the study of the convergence of $(V_{\theta^k})_k$.

\begin{cor}\label{corosupinf}  For all $T_0\geq 0$ and $y_0$ in $\RR^d$,
$$ \sup_{k\geq 1}\inf_{t\leq T_0} V_{\T_t\sharp\theta^k}(y_0)\geq \limsup_{k}V_{\theta^k}(y_0)\geq \liminf_{k}V_{\theta^k}(y_0)\geq \sup_{k\geq 1}\inf_{t\geq 0} V_{\T_t\sharp\theta^k}(y_0)$$
\end{cor}
\begin{rem} Corollary \ref{corosupinf} states that the uniform convergence of "$\sup_{k\geq 1}\inf_{t\leq T_0} V_{\T_t\sharp\theta^k}$" to "$\sup_{k\geq 1}\inf_{t\geq 0} V_{\T_t\sharp\theta^k}$" as $T_0$ tends to infinity implies the uniform convergence of $(V_{\theta^k})_k$ as $k$ tends to infinity. Moreover, according to Corollary \ref{coro:V*}, in case of uniform convergence, the limit function is $V^*$.
\end{rem}
\noindent For any states $y$ and $\overline{y}$ in $\RR^d$, let us define $\tilde{d}(y,\overline{y})=\sup_{k\geq 1} |V_{\theta^k}(y)-V_{\theta^k}(\overline{y})|$. The space $(\RR^d,\tilde{d})$ is now a \textit{pseudometric} space (may not be Hausdorff).

\noindent The following is similar to the proof of Theorem 2.5 in \cite{Renault_2014}, and is also similar to the proof of Theorem 3.10 in \cite{Renault_2011}. We rewrite it here for sake of completeness. Roughly speaking, we shall use the total boundedness of the space $\big(\{V_{\theta^k}\},||\cdot||_\infty\big)$ so as to deduce that the state space $(\RR^d, \tilde{d})$ is totally bounded for the pseudometric metric $\tilde{d}$. This allows us to prove the convergence for $\tilde{d}$ of the reachable set $R^T$ to $R$ in bounded time. We are then able to prove the uniform convergence of "$\sup_{k\geq 1}\inf_{t\leq T_0}V_{\T_t\sharp \theta}$" to "$\sup_{k\geq 1}\inf_{t\geq 0}V_{\T_t\sharp \theta}$" as $T_0$ tends to infinity..

\noindent{\bf{Proof for Theorem \ref{thm:main}, Parts (ii)-(iii)}}.

We first prove Part (iii). One direction is clear: the uniform convergence of $(V_{\theta^k})$ implies the totally boundedness of the space $(\{V_{\theta^k}\},||\cdot||_\infty)$.

Let us prove the converse. Suppose that $(\{V_{\theta^k}\}, ||\cdot||_\infty)$ is totally bounded, so fixing any $\varepsilon>0$, there exists a finite set of indices $I$ such that for all $k\geq 1$, there exists $i\in I$ satisfying $$||V_{\theta^k}-V_{\theta^i}||_\infty\leq \varepsilon/3.$$
$\{\big(V_{\theta^i}(y)\big),y\in \RR^d\}$ is a subset of the compact metric space $\big([0,1]^I,\|\cdot\|_\infty\big)$, thus it is itself totally bounded, so there exists a finite subset $X$ of $\RR^d$ such that $$\forall y\in \RR^d, \exists x\in X, \forall i\in I, |V_{\theta^i}(y)-V_{\theta^i}(x)|\leq\varepsilon/3. $$
We have obtained that for each $\varepsilon>0$, there exists a finite subset $X$ of $\RR^d$ such that for every $y\in\RR^d$, there is $x\in X$ satisfying: for any $k\geq 1$ there is some $i\in I$ with
 $$\big|V_{\theta^k}(y)-V_{\theta^k}(x)\big|\leq \big|V_{\theta^k}(y)-V_{\theta^i}(y)\big|+\big|V_{\theta^i}(y)-V_{\theta^i}(x)\big|+\big|V_{\theta^i}(x)-V_{\theta^k}(x)\big|\leq \varepsilon,$$
thus $\tilde{d}(y,x)\leq \varepsilon$. This implies that the pseudometric space $(\RR^d,\tilde{d})$ is itself totally bounded.

\noindent Fix now $y_0$ in $\RR^d$. It is by definition that $$\text{ for all }  \ T, S\in\RR_+ \text{ with } S\geq T, \text{ we have } R^{T}(y_0)\subset R^{S}(y_0)\subset R(y_0),$$ and
 $$\forall \bar{y}\in R(y_0), \ \exists \bar{T}>0 \text{ with } \bar{y}\in R^{\bar{T}}(y_0).$$

From the totally boundedness of $(\RR^d,\tilde{d})$, we show that $R^T$ converges to $R$ in the following sense
\begin{equation} \label{eq:R^T}
\forall \varepsilon>0, \ \exists T\geq 0: \ \forall \overline{y}\in R(y),\  \exists \tilde{y}\in R^T(y), \ \tilde{d}(\overline{y}, \tilde{y})\leq \varepsilon.
\end{equation}

Indeed, let us first take $\{y_\ell\}$ a finite $\varepsilon$-cover of $R(y_0)$ for $\tilde{d}$. For each $y_\ell$, put $T_\ell>0$ with $y_\ell\in R^{T_\ell}(y_0)$. We then take $T=\max T_\ell$. Now for any $\bar{y}\in R(y_0)$, there is some $y_\ell$ with $\tilde{d}(\bar{y}, y_\ell)\leq \varepsilon$. Moreover, $y_\ell\in R^{T_\ell}(y_0)\subset R^T(y_0)$. This proves (\ref{eq:R^T}).

\noindent By Corollary \ref{corosupinf}, for all $T\geq 0$ \big(\ using $\inf_{\overline{y}\in R^T(y_0)} V_{\theta^k}(\overline{y})=\inf_{t\leq T} V_{\T_t\sharp \theta^k}(y_0)$, cf. (\ref{eq:V+t+shift})\big):
$$ \sup_{k\geq 1}\inf_{\overline{y}\in R^T(y_0)} V_{\theta^k}(\overline{y})\geq \limsup_{k}V_{\theta^k}(y_0)\geq \liminf_{k}V_{\theta^k}(y_0)\geq \sup_{k\geq 1}\inf_{\overline{y}\in R(y_0)} V_{\theta^k}(\overline{y}).$$
\noindent Consider $k\geq 1$ and $T\geq 0$ given by assertion (\ref{eq:R^T}) for the fixed $\varepsilon>0$. Let $\overline{y}\in R(y_0)$ be such that $V_{\theta^k}(\overline{y}_0)\leq \inf_{\overline{y}\in R(y)}V_{\theta^k}(\overline{y})+\varepsilon$, and then $\tilde{y}$ in $R^T(y_0)$ be such that $\tilde{d}(\overline{y},\tilde{y})\leq \varepsilon$. Since $V_{\theta^k}$ is clearly 1-Lipschitz for $\tilde{d}$, we obtain $V_{\theta^k}(\tilde{y})\leq \inf_{\overline{y}\in R(y_0)}V_{\theta^k}(\overline{y})+2\varepsilon$. Consequently, $\inf_{\overline{y}\in R^T(y_0)}V_{\theta^k}(\overline{y})\leq \inf_{\overline{y}\in R(y_0)} V_{\theta^k}(\overline{y})+2\varepsilon$ for all $k$, so
$$\sup_{k\geq 1}\inf_{\overline{y}\in R^T(y_0)}V_{\theta^k}(\overline{y})\leq \sup_{k\geq 1}\inf_{\overline{y}\in R(y_0)}V_{\theta^k}(\overline{y}) +2\varepsilon.$$
One obtains that $\limsup_{k\geq 1} V_{\theta^k}(y_0)\leq \liminf_{k\geq 1} V_{\theta^k}(y_0)+2\varepsilon$, and so $\big(V_{\theta^k}(y_0)\big)_k$ converges. Since $(\RR^d,\tilde{d})$ is totally bounded and all $V_{\theta^k}$ are 1-Liptschitz, the convergence is uniform.

\noindent Next, Part (ii) can be deduced from the proof of Part (iii). Let $(\theta^{\varphi(k)})$ be any subsequence of $(\theta^k)$ that converges uniformly to some function $V$. This implies that $(\{V_{\theta^{\varphi(k)}}\},||\cdot||_\infty)$ is totally bounded. As we have shown in the proof of Part (iii) that if $(\{V_{\theta^{\varphi(k)}\}},||\cdot||_\infty)$ is totally bounded, $(V_{\theta^{\varphi(k)}})$ converges uniformly to $V=V^*$, which implies Part (ii) that $V^*$ is the unique accumulation point (for the uniform convergence) of the sequence $(V_{\theta^k})_k$.
$\hfill \Box$

\section{Discussion on a weaker long-term condition}\label{sec:LTC'}

Below is a weaker form of the long-term condition (LTC):

\noindent \textbf{Long-term condition' (LTC')} \ A sequence of evaluations $(\theta^k)_{k\geq 1}$ satisfies the LTC' if:
\begin{align} \label{eq:LTC'}
\forall s>0, \ \ TV_s(\theta^k)\xrightarrow[k \to \infty]{} \ 0.
\end{align}

It is unclear whether the LTC' is strictly weaker than the LTC or not. One might want to construct an example of $(\theta^k)_k$ such that $TV_s(\theta^k) \xrightarrow [k\to\infty]  {} \ 0  $ for all $s>0$ while $\overline{TV}_{s_0}(\theta^k)\xrightarrow [k\to \infty]  {} \ \alpha>0$ for some $s_0>0$ and $\alpha>0$. The following example shows that this is possible if we consider only $s$ being rational numbers. In general, the question is still open.

\begin{exa}\label{exa:LTC'} Given a positive integer $k$, consider the density $\theta^k$ with support included in $[0,k]$ by   dividing  $[0,k]$  in $k^2$ consecutive small intervals of length $1/k$, and $\theta^k$ is uniform over the union of all small odd intervals... and puts no weight on even small intervals. Define the support $$S_k=\bigcup_{l \in \NN, l\leq \frac{k^2-1}{2}, } \left[\frac{2l}{k}, \frac{2l+1}{k}\right).$$
$\theta^k$ has density:$$f_k(x)=\frac{2}{k}\;  {\mathds{1}}_{x \in S_k}= \frac{2}{k}\;  {\mathds{1}}_{x\in [0,k], E(kx) \in 2\NN}$$
(where $2 \NN$ is the set of even numbers in $\NN$, $E(x)$ is the integer part of $x$).

For each $k$, we have (consider $s=1/k$):
 $$\sup_{0\leq s\leq 1} \int_{x\geq 0} \big| f_k(x+s)-f_k(x)\big| \deriv x \geq  2-1/k$$
Consider now only $k$ of the form $n!$, and we define the density $g_n=f_{n!}$ for each $n$ in $\NN$. For all $x\geq 0$,
$$g_n(x+s)-g_n(x)=\frac{2}{n!} \left({\mathds{1}}_{E(n!(x+s))\in 2 \NN, x+s\leq n!}- {\mathds{1}}_{E(n! x)\in 2 \NN, x\leq n!}\right).$$
Assume $s$ is a rational number. Then for $n$ large enough, $n!s$ is an even  integer, so for all $x$ such that $0\leq x \leq n!-s$, we have $g_n(x+s)-g_n(x)=0$. Consequently,  $$ \int_{x\geq 0} \big| g_n(x+s)-g_n(x)\big| \deriv x \xrightarrow [n \to \infty]{} \ 0.$$
\end{exa}

\section{Appendix}\label{sec:appendix}

\noindent \textbf{Proof for Lemma \ref{lem:variation}}: The following computation of $f'_{\theta}(t)$ is straightforward:
$$ \ \forall t>0,\ \  f'_\theta(t)= \frac{1}{\sigma \sqrt{2\pi}}\Big[\exp\left(-\frac{1}{2}\Big(\frac{t-m}{\sigma}\Big)^2\right)\frac{m-t}{\sigma^2}- \exp\left(-\frac{1}{2}\Big(\frac{t+m}{\sigma}\Big)^2\right)\frac{m+t}{\sigma^2}\Big],$$
thus
$$ f'_\theta(t) >0  \ \ ( resp. <0)
\Longleftrightarrow  \ (m-t)\exp\left(-\frac{1}{2}\Big(\frac{t-m}{\sigma}\Big)^2\right) - (m+t) \exp\left(-\frac{1}{2}\Big(\frac{t+m}{\sigma}\Big)^2\right) >0\  \ ( resp. <0).$$

As a consequence, one obtains that $f'_\theta(t)< 0,\ \forall t\geq m$. Now we look at $t\in(0,m)$. Denote $H(t)=_{def}\exp\left(\frac{2mt}{\sigma^2}\right)-\frac{m+t}{m-t}$, which enables us to write: $$f'_\theta(t)>0 \ ( resp. <0) \Longleftrightarrow  H(t)>0  \ ( resp. <0),  \ \ \forall t\in(0,m).$$
From the above analysis, we deduce that the proof of the lemma is reduced to the proof for

\noindent \textbf{Claim}\  \textit{There is some $t^*\in[0,m)$ such that $H(t)<0$ for $t\in(0,t^*)$ and $H(t)>0$ for $t\in (t^*, m)$. Moreover, such $t^*$ satisfies $(t^*)^2\geq m^2-\sigma^2$.}

In order to prove the claim, we compute
\begin{itemize}
\item the values at the end point: $H(0)=0$ and $\lim_{t\to m^-}H(t)=-\infty$;
\item the first-order derivative at any $t\in[0,m)$:
\begin{align} \label{eq:H'}
H'(t)= \exp\left(\frac{2mt}{\sigma^2}\right)\frac{2m}{\sigma^2}-\frac{2m}{(m-t)^2}
\end{align}
\item at any rest point $t^e\in[0,m)$ ($i.e.$, $H(t^e)=0$):
\begin{align*}
\exp\left(\frac{2mt^e}{\sigma^2}\right)=\frac{m+t^e}{m-t^e},
\end{align*}
which is substituted back into (\ref{eq:H'}), to yield
\begin{align}\label{eq:sharp}
H'(t^e)>0 \  \big(\ resp.\ H'(t^e)<0 \ \big) \  \Longleftrightarrow  \ (t^e)^2 <m^2-\sigma^2  \ \big(\ resp.\  (t^e)^2>m^2-\sigma^2 \ \big).
\end{align}
\end{itemize}
Next, it is easy for us to prove the following result:

\textit{Let $t^e_1\in[0,m)$ be a rest point for $H(\cdot)$, and suppose that $t^e_2\in(t^e_1,m)$ is the smallest rest point after $t^e_1$. Then $H'(t^e_1)H'(t^e_2)\leq 0$ and if $H'(t^e_1)\leq 0$,  such $t^e_2$ does not exist.}

Indeed,  $H'(t^e_1)H'(t^e_2)\leq 0$ can be derived from the continuity of $H(\cdot)$; suppose that $H'(t^e_1)\leq 0$, we have from (\ref{eq:sharp}) that $(t^e_1)^2\geq m^2-\sigma^2$ and $H'(t^e_2)\geq 0$, thus $(t^e_2)^2 \leq m^2-\sigma^2$. However, this leads to a controdiction to $t^e_2>t^e_1$, so $t^e_2$ does not exist whenever $H'(t^e_1)\leq 0$.

Finally, remark that $H(0)=0$, thus $t=0$ is a rest point. We discuss the following two cases:

\underline{Case 1. $m^2-\sigma^2\leq 0$, thus $H'(0)\leq 0$}.

This implies that no rest point exists after $0$. Since $\lim_{t\to m^-}H(t)=-\infty$, we deduce that $H(t)<0, \ \forall t\in(0,m)$. The claim is proved for $t^*=0$.

\underline{Case 2. $m^2-\sigma^2> 0$, thus $H'(0)>0$}.

$\lim_{t\to m^-}H(t)=-\infty$ implies that some rest point exists in $(0,m)$. Take $t^e$ the closest to $0$, implying that $H(t)> 0, \ \forall t\in(0,t^e)$. Further, we obtain that $H'(t^e)\leq 0$ by the continuity of $H(\cdot)$. Again, there exists no other rest point after $t^e$. Since $\lim_{t\to m^-}H(t)=-\infty$, we deduce that $H(t)<0, \ \forall t\in(t^e,m)$. The claim is proved for $t^*=t^e$.

To conclude, we see that in both cases such $t^*$ exists and satisfies $(t^*)^2\geq m^2-\sigma^2$, thus the claim is proved. This finishes our proof for the lemma. $\hfill\Box$

\bibliographystyle{plain}
\bibliography{Control_Li-Quincampoix-Renault}

\noindent \textit{Email address}: xxleewhu@gmail.com  \\
\textit{Email address}: Marc.Quincampoix@univ-brest.fr\\
\textit{Email address}: jerome.renault@tse-fr.eu
\end{document}